\documentclass[11pt]{amsart}
\usepackage{amsfonts,amssymb,amsthm,amsmath,amsxtra,amscd,verbatim,eucal,pdfsync}
\usepackage[all]{xy}
\usepackage[dvips]{graphics}

\theoremstyle{plain}
\newtheorem{thm}{Theorem}[section]
\newtheorem{lem}[thm]{Lemma}
\newtheorem{prop}[thm]{Proposition}

\newtheorem{theoremalpha}{Theorem}
\newtheorem{corollaryalpha}[theoremalpha]{Corollary}
\newtheorem{propositionalpha}[theoremalpha]{Proposition}

\theoremstyle{definition}

\newtheorem{eg}[thm]{Example}

\theoremstyle{remark}
\newtheorem{rmk}[thm]{Remark}

\newtheorem{nothing}[thm]{}

\def\Z{{\mathbb Z}}

\def\F{{\mathbb F}}
\def\C{{\mathbb C}}

\def\R{{\mathbb R}}
\def\Q{{\mathbb Q}}
\def\P{{\mathbb P}}

\def\O{\mathcal{O}}

\def\aa{\mathfrak{a}}
\def\bb{\mathfrak{b}}

\def\mm{\mathfrak{m}}

\def\a{\alpha}

\def\d{\delta}

\def\ep{\varepsilon}

\def\l{\lambda}
\def\n{\nu}
\def\m{\mu}

\def\r{\rho}

\def\.{\cdot}
\def\^{\widehat}
\def\~{\widetilde}

\def\({\left(}
\def\){\right)}

\DeclareMathOperator{\im}  {im}

\DeclareMathOperator{\Bl} {Bl}

\DeclareMathOperator{\Div} {Div}

    \newcommand{\eps}{\varepsilon}
\newcommand{\BBB}{\mathbf{B}}

\newcommand{\HH}[3]{H^{{#1}} \big( {#2} , {#3}
\big) }
\newcommand{\hh}[3]{h^{{#1}} \big( {#2} , {#3}
\big) }
\newcommand{\HHnospace}[2]{H^{{#1}} \big( {#2}   \big) }
\newcommand{\hhnospace}[2]{h^{{#1}} \big( {#2}   \big) }
\newcommand{\andd}{ \ \ \text{ and } \ }

\newcommand{\fall}{ \ \ \text{ for all } \ }
\newcommand{\hhat}[2]{\widehat{h}^{{#1}}\big( {#2}
\big)}
\newcommand{\Linser}[1]{| \mspace{1.5mu} {#1}
\mspace{1.5mu} |}
\newcommand{\linser}[1]{\Linser{  {#1}  }}

\newcommand{\bs}[1]{\bb \big(\mspace{1.1mu}
\linser{{#1}} \mspace{1.1mu} \big)}
\newcommand{\oversetoplus}[1]{\overset{#1}{\oplus}}
\newcommand{\noi}{\noindent}
\newcommand{\lra}{\longrightarrow}
\def\OO{\mathcal{O}}
\newcommand{\MI}[1]{\mathcal{I} \big( {#1} \big) }

\newcommand{\vol}{\textnormal{vol}}

\newcommand{\deq}{\ensuremath{
\stackrel{\textrm{def}}{=}}}

\begin{document}
\title{Higher cohomology of divisors on a projective variety}

\author{Tommaso de Fernex}
\address{Department of Mathematics, University of Utah, 155 South 1400
East, Salt Lake City, UT 48112-0090}
\curraddr{Institute for Advanced Study, School of Mathematics, 1 Einstein
Drive, Princeton, NJ 08540}
\email{{\tt defernex@math.utah.edu}}

\author{Alex K\"uronya}
\address{Mathematical Institute, Budapest University of Technology and
Economics, P.O. Box 91, Budapest, H-1152 Hungary}
\address{Universit\"at Duisburg-Essen, Campus Essen, FB 6 Mathematik,
D-45117 Essen, Germany}
\email{{\tt alex.kuronya@math.bme.hu}}

\author{Robert Lazarsfeld}
\address{Department of Mathematics, University of Michigan, Ann Arbor, MI 48109, USA}
\email{{\tt rlaz@umich.edu}}

\thanks{Research of the first author partially supported by NSF grants
DMS 0111298 and DMS 0548325,
and by MIUR National Research Project
``Geometry on Algebraic Varieties" (Cofin 2004).
Research of the second author was partially supported by the Leibniz program
of the Deutsche Forschungsgemeinschaft and the OTKA grant T042481 of the Hungarian
Academy of Sciences. Research of the third author partially
supported by NSF grant DMS 0139713.}

\maketitle

\section*{Introduction}

The purpose of this paper is to study the growth of higher cohomology of
line bundles on a projective variety.

Let $X$ be an irreducible complex projective variety of dimension $d \ge
1$,  and let $L$ be a Cartier divisor  on $X$.  It is elementary and
well-known that the dimensions of the cohomomology groups
$\HH{i}{X}{\OO_X(mL)}$ grow at most like $m^d$, i.e.
$$
\dim_\C \HHnospace{i}{X,\O_X(mL)} \ = \ O(m^d)
\fall \ i \ge 0
$$ (cf. \cite[1.2.33]{PAG}). It is natural to ask when one of these
actually has maximal growth, i.e. when
$\hh{i}{X}{\O_X(mL)}   \ge   Cm^d$ for some positive constant $C > 0$ and
arbitrarily large $m$.
 For $i = 0$ this happens by definition exactly when
$L$ is \textit{big}, and the geometry of big classes is fairly well
understood. Here we focus on the question of when one or more of the
higher cohomology groups grows maximally.

  If $L$ is ample, or merely nef,   then of course
$\hh{i}{X}{\OO_X(mL)} = o(m^{d})$ for $i > 0$. In general the converse is
false: for instance, it can happen that $\HH{i}{X}{\OO_X(mL)} = 0$ for
$m > 0$ and all $i$ even if $L$ is not nef (or, for that matter,
pseudoeffective).
     However our main result shows that if one considers also small
perturbations of the divisor in question, then in fact the maximal
growth  of higher cohomology characterizes non-ample divisors:
    \begin{theoremalpha} \label{Max.Higher.Cohom}   Fix any very  ample
divisor $A$ on $X$.  If
$L$ is not ample, then  for     sufficiently small rational numbers
$t > 0$, at least one of the higher cohomology groups of  suitable
multiples of  $L - tA$ has maximal growth. More precisely, there is an
index $i > 0$ such that for any sufficiently small $t > 0,$
\[
\dim_\C \HHnospace{i}{X,\O_X(m(L-tA))} \ \ge \ C \cdot m^d  \] for some
constant $C = C(L,A,t) > 0 $ and  arbitrarily large values of $m$
clearing the denominator of $t$.
\end{theoremalpha}
\noi In other words,  a divisor $L$ is ample if and only if
\[  \hh{i}{X}{\OO_X(m(L - tA))} \ =  \ o(m^d) \ \  \text{when  } i > 0\]
for all small $t$ and suitably divisible $m$.   We remark that the
essential content of the theorem is the maximal growth of higher
cohomology when $L$ is big but not ample.

One can get a more picturesque statement by introducing asymptotic
invariants of line bundles. As above let $X$ be an irreducible complex
projective variety of dimension $d$, and let $L$ be a Cartier divisor on
$X$. Define
$$
\^h^i(X,L) \ \deq \ \limsup_{m \to \infty}
\frac{\dim_\C \HHnospace{i}{X,\O_X(mL)}}{m^d/d!}.
$$  The definition extends in the natural manner to
$\Q$-divisors. When
$i = 0$ this is the
\textit{volume}
$\vol_X(L)$ of
$L$, which has been the focus of considerable attention in recent years
(\cite{Fuj}, \cite{DEL},
\cite{Boucksom}, \cite{PAG}). The higher cohomology functions were
introduced and studied by the second author in \cite{Kur}. It was
established there that
$\^h^i(L)$ depends only on the numerical equivalence class of a
$\Q$-divisor $L$, and that it uniquely determines a continuous function
\[
\^h^i =  \hhat{i}{X,\;\;\;} \colon N^1(X)_\R \lra \R\] on the real
N\'eron-Severi space of $X$. This generalizes the corresponding results
for $\vol_X$ proved by the third author in \cite[2.2.C]{PAG}. When $X$ is
a toric variety, these functions were studied in \cite{HKP}.  We refer to
\cite{ELMNP} for a survey of the circle of ideas surrounding asymptotic
invariants.

Theorem \ref{Max.Higher.Cohom} then implies:
\begin{corollaryalpha}
\label{Characterize.Ample.Classes.Cor.Intro}
 A class $\xi_0 \in N^1(X)_{\R}$ is ample if and
only if
\[
\^h^i(\xi ) \ = \ 0 \] for all $i > 0$ and all $\xi  \in N^1(X)_{\R}$ in
a small neighborhood of $\xi_0$.
\end{corollaryalpha}
\noi One can see this as an asymptotic analogue of Serre's criterion for
amplitude. In the toric case, a somewhat stronger statement appears in
\cite{HKP}.

The proof of the theorem combines some algebraic constructions from
\cite{Kur} with  geometric facts about big line bundles that fail to be
ample. Specifically, choose a collection of very general divisors $E_{1},
\ldots, E_p \in \linser{A}$, and assume for a contradiction that $L$ is
not nef, but that
\[  \^h^{i}\big(X, L - tA\big)   \ =  \ 0 \] for $i >0$ and small
positive $t$. The first point is to show that this vanishing descends to
the divisors $E_{\alpha}$. By induction it follows each of the
restrictions $L|_{E_\alpha}$ is ample, and we can
 assume furthermore that they are all very ample. Now consider the complex
\begin{equation}
\HHnospace{0}{\OO_X(mL)}
\overset{v}{\lra}\oversetoplus{p}\,
\HHnospace{0}{\OO_{E_{\alpha}}(mL)}\overset{u}{\lra}
\oversetoplus{\binom{p}{2}}\,\HHnospace{0}{\OO_{E_\alpha
\cap E_\beta}(mL)} ,
\tag{*}
\end{equation}  the maps $u$ and $v$ being determined  by restriction.
The cohomology of this injects into
$\HH{1}{X}{\OO_X(mL - pA)}$. On the other hand, using the assumption that
$L$ fails to be  nef, we show that one can arrange things so that there
is a non-trivial ideal $\aa \subseteq \OO_X$, vanishing on a set of
dimension
$\ge 1$, such that the base ideal of $\linser{mL}$ grows like $\aa^m$,
i.e.
\[
\bs{mL} \ \subseteq \ \aa^m.
\] Therefore the image of $v$ is contained in the subgroup
\[
\oversetoplus{p}\,\HHnospace{0}{\OO_{E_{\alpha}}(mL)
\otimes\aa^m} \ \subseteq \
\oversetoplus{p}\,\HHnospace{0}{\OO_{E_{\alpha}}(mL)},
\] and this leads to a lower bound on the dimension of $\ker(u)/\im(v)$.
In fact,  a dimension count shows that if $p \sim  m\delta$ for $0
<\delta\ll 1$, then $\hhnospace{1}{mL - pA}$ will grow like $m^d$, which
produces the required cohomology.

This argument suggests that there is a relation between the amplitude of
restrictions of $L$ and its asymptotic higher cohomologies, and we
explore this connection in Section 3. Fixing as above a very ample
divisor $A$, define $a(L,A)$ to be the least integer $k$ such that
$L\mid _{E_1 \cap \ldots \cap E_k} $ is ample for $k$ very general
divisors $E_1, \ldots, E_k \in \linser{A}$. So for example, $a(L,A) \le
\dim \BBB_+(L)$, where $\BBB_+(L)$ is the augmented stable base locus of
$L$ in the sense of \cite{AIBL}, (cf. \cite[Section 10.3]{PAG}).\footnote{By definition,
$\BBB_+(L)$ is the stable base-locus of the $\Q$-divisor $L - \eps A$ for
any small $\eps > 0$ and ample $A$.} We prove
\begin{propositionalpha} If $i > a(L,A)$ then $\^h^{i}(\xi) = 0$ for all
$\xi$ in a neighborhood of $[L]$ in $N^1(X)_{\R}$.
\end{propositionalpha}

\noi We also give examples to show that the converse can fail, i.e. that
it can happen that $\^h^{i}(L - tA) = 0$ for small $t > 0$ when $i =
a(L,A)$. It remains an interesting open question whether one can predict
from geometric data the largest value of $i$ for which
$\^h^{i}(L - tA) \ne  0$ for small $t>0$.

The paper is organized as follows.
We start in Section 1 with a lemma concerning base ideals of big linear series.
The main result appears in Section 2, and finally in Section 3 we study restrictions
and asymptotic higher cohomology.

The first author would like to thank Eckart Viehweg for support during his
visit at the Universit\"at Duisburg-Essen. We are grateful to Lawrence Ein
and Sam Payne for useful discussions.

\setcounter{section}{-1}

\section{Conventions and  Background}

\begin{nothing} We work throughout over the complex numbers. A
\textit{variety} is a reduced and irreducible scheme, and we always deal
with closed points.
\end{nothing}

\begin{nothing} We follow the conventions of \cite[Chapter 1]{PAG} concerning
divisors on a projective variety $X$. Thus a \textit{divisor} on $X$ means
a Cartier divisor. A $\Q$- or $\R$-divisor   indicates an element of
\[  \Div_{\Q}(X) \ = \ \Div(X) \otimes \Q \ \   \text{or} \ \
\Div_{\R}(X) =
\Div(X) \otimes \R. \] $N^1(X)$ is the N\'eron-Severi group of numerical
equivalence classes of divisors, while $N^1(X)_{\Q}$ and  $N^1(X)_{\R}$
denote the corresponding groups for $\Q$- and $\R$-divisors.
\end{nothing}

\begin{nothing}
Given a projective variety $X$ of dimension $d$, and a divisor $L$ on $X$,
we set
\begin{equation}
\^h^i(X,L) \ = \ \limsup_{m \to \infty}
\frac{\dim_\C \HHnospace{i}{X,\O_X(mL)}}{m^d/d!}.
\end{equation}
It is established in \cite[Proposition 5.15, Theorem 5.1]{Kur} that this
depends only on the numerical equivalence class of $L$, and that it
satisfies the homogeneity
\begin{equation} \label{Homogenity.Cohom.eqn}
\^h^i(pL) \ = \ p^d \cdot \^h^i(L).
\end{equation}
This allows one to define $\^h^i(L)$ for an arbitrary $\Q$-divisor by
clearing denominators, giving rise to a function
$ \^h^i : N^1(X)_{\Q} \lra \R.$ The main result of \cite{Kur} is that
this extends uniquely to a continuous function
\[ \^h^i : N^1(X)_{\R} \lra \R \]
satisfying the same homogeneity property (\cite[Theorem 5.1]{Kur}). When
$L$ is an
$\R$-divisor, we typically write $\^h^i(L)$ to denote the value of this
function on the numerical equivalence class of $L$. The $\widehat{h}^i$ are
called asymptotic cohomological
functions in \cite{Kur}, although we occasionally use some slight variants
of this terminology.
 Observe finally that
the homogeneity  \eqref{Homogenity.Cohom.eqn} implies that if $L$ is a
Cartier divisor, then
\begin{equation}
\^h^i(X,L) \ = \ \frac{\^h^i(X,pL) }{p^d} \ = \ \limsup_{m \to \infty}
\frac{\dim_\C \HHnospace{i}{X,\O_X(pmL)}}{(pm)^d/d!}.
\end{equation}
for every fixed integer $p > 0$. The analogous statement holds when $L$
is a $\Q$-divisor, provided that $p$ is sufficiently divisible to clear
the denominators of $L$.

\end{nothing}

\section{A lemma on base-loci}

This brief section is devoted to the proof of a useful fact concerning
the base-ideals of linear systems on a normal variety.
The result in question asserts  that the base-ideals
associated to multiples of a  divisor which is not nef grow at least
like  powers of the ideal of a curve. On a smooth variety this fact,
which is a small elaboration of \cite[Corollary~11.2.13]{PAG},
is well-known to experts. The main point here is to extend the statement
to normal varieties.\footnote{We thank Mircea Musta{\c t}{\v a}
for pointing out that in a first version of this paper,
this reduction contained an error.}

\begin{prop}\label{lemma:base-loci} Let $D$ be a divisor on a normal
projective variety $V$, and denote by \[ \bs{\ell D} \subseteq \OO_V\]
the base-ideal of the indicated linear series. Assume that  $D$ is not
nef. Then there exist positive integers $q$ and $ c$, and an ideal sheaf
$\aa \subseteq \OO_V$ vanishing on a set of dimension $\ge 1$, such that
\[
\bs{mqD} \ \subseteq \ \aa^{m-c}
\] for all $m > c$.
\end{prop}

\begin{proof} The assertion being trivial otherwise, we may suppose that
$D$ has non-negative Kodaira-Iitaka dimension $\kappa(V, D) \ge 0$.
Since $D$ is not nef, there exists an irreducible curve $Z \subset V$
intersecting $D$ negatively. We will show that one can take  $\aa = I_Z \subset \O_V$ to
be the ideal sheaf of $Z$.

Let $\m : V' \to V$ be a resolution of singularities,
and set $D' = \m^*D$. We may suppose that
  $$
\aa \cdot\O_{V'} \ = \ \O_{V'}(-E),
$$
where $E$ is an effective divisor on $V'$ with simple normal crossing support.
Note that the projection formula implies that
 $D'$ intersects negatively every curve $C \subset V'$ that dominates $Z$.

Given $q \ge 1$, consider the asymptotic multiplier ideal sheaf
\[ \MI{V^\prime, \Vert qD' \Vert} \ =\
\MI{ \Vert qD' \Vert} \ \subseteq \ \OO_{V^\prime}. \]
Since $D'$ has non-negative Kodaira-Iitaka dimension,
there exists a fixed divisor $A$ on $V'$ such that
\[
\OO_{V'}(qD' + A) \, \otimes \, \MI{\Vert qD' \Vert}
\]
is globally generated for every $q$;
this follows from \cite[Corollary~11.2.13]{PAG}
by taking $A = K_{V'} + (\dim V' + 1)H$, where $H$
is any very ample divisor on $V'$.
If $C \subset V'$ is an irreducible curve that
is not contained in the zero locus of
$\MI{\Vert qD' \Vert}$, then one has
\[
\big((qD' + A) \cdot C \big)\ \ge\ 0.
\]
We conclude that if $C \subseteq V^\prime$ is any curve dominating $Z$,
then the ideal $\MI{\Vert qD' \Vert}$
must vanish along  $C$ for  $q > -(A \cdot C)/(D^\prime \cdot C)$.

\begin{lem}\label{lem:vanishing-of-mult-ideal}
There exists $q > 0$ such that $\MI{\Vert qD' \Vert} \subseteq \O_{V'}(-E)$.
\end{lem}

\begin{proof}
We observe that it is enough to prove that
\begin{equation}\label{eq:E-reduced}
\MI{\Vert kD' \Vert} \subseteq \O_{V'}(-E_{\rm red}) \ \ \text{for some} \ k > 0.
\end{equation}
Indeed, assuming this, if $k$ is as
in \eqref{eq:E-reduced} and $a$ is the largest
coefficient appearing in the divisor $E$, then by the subadditivity
theorem for multiplier ideals (\cite[Corollary~11.2.4]{PAG})  we have
$$
\MI{\Vert akD' \Vert}
\ \subseteq \ \MI{\Vert kD' \Vert}^a \ \subseteq \ \O_{V'}(-a E_{\rm red})
\ \subseteq \ \O_{V'}(-E).
$$
Hence we can take $q = ak$.

 Turning to \eqref{eq:E-reduced},  note to begin with  that for
 some $k_0$, the ideal sheaf $\MI{\Vert k_0 D' \Vert}$
vanishes along all irreducible components of $E$ that
dominate $Z$. In fact, suppose that $F$ is such a component,
and let $C \subset F$ be a general complete intersection curve
that surjects onto $Z$. Then as we have seen, $\MI{\Vert qD' \Vert}$
vanishes along $C$ for $q \gg 0$. So this ideal must in fact vanish on all of $F$.

Now since $V$ is normal, $E = \mu^{-1}(Z)$ is connected.
Therefore, in order to prove \eqref{eq:E-reduced},
it suffices to prove the following
claim: If $E_i$ is an irreducible component of $E$
mapping to a point of $Z$,  and if $E_i$ meets another irreducible component $E_j$
of $E$ along which $\MI{\Vert k_j D' \Vert}$ is known to vanish for
some $k_j > 0$, then there
is a $k_i \ge k_j$ such that $\MI{\Vert k_i D' \Vert}$ vanishes along $E_i$.
(Note that then $\MI{\Vert k_i D' \Vert}$ also vanishes along $E_j$
since $\MI{\Vert \ell  D' \Vert} \subseteq \MI{\Vert m  D' \Vert}$ whenever $\ell > m$.)

So, let $E_i$ and $E_j$ be as in this scenario.
Since $\MI{\Vert k_j D' \Vert}$ vanishes on $E_j$, we have
$$
\MI{\Vert mk_j D' \Vert} \ \subseteq \ \MI{\Vert k_j D' \Vert}^m \
\subseteq \ \O_{V'}(-mE_j) \fall m > 0
$$
thanks again to the subadditivity theorem.
Now take a general complete intersection curve $C \subseteq E_i$ that meets $E_j$ in
at least one point $P$, which we may assume to be a smooth point of $C$.
As above, we know that
$$
\O_{V'}(mk_j D'+A) \ \otimes \ \MI{\Vert mk_j D' \Vert}
$$
is globally generated for $m > 0$.
On the other hand, $\OO_C(D^\prime) = \OO_C$ since $C$ maps to a point in $V$,
and hence $$
\Big(\O_{V'}(mk_j D'+A) \ \otimes \ \MI{\Vert mk_j D' \Vert}\Big) \ \cdot \ \O_C
\ \subseteq  \ \O_C(A|_C - mP).$$
Therefore, assuming that $\MI{\Vert mk_j D' \Vert}$ does not vanish along $C$,
it follows that $(A \cdot C) \ge m$.
  Thus if  $k_i \ge k_j\cdot\big(  (A \cdot C)+ 1 \big)$,
then in fact $\MI{\Vert k_i D' \Vert}$ must vanish along $E_i$.
This completes the proof of the lemma.
\end{proof}

We now complete the proof of the Proposition. Note that
$\HH{0}{V}{\OO_V(\ell D)}=\HH{0}{V^\prime}{\OO_{V^\prime}(\ell D^\prime)}$
by normality, and in particular
$\bs{\ell D^\prime} = \bs{\ell D} \cdot \OO_{V^\prime}$.
On the other hand, we have
\[
\bs{mqD'} \ \subseteq \  \MI{\Vert mqD' \Vert}
\ \subseteq \ \MI{\Vert qD' \Vert}^m \ \subseteq \ \O_{V'}(-mE)
\]
by \cite[Corollary~11.2.4]{PAG} and Lemma~\ref{lem:vanishing-of-mult-ideal}.
Therefore
\[
\bs{mq D} \ \subseteq \ \mu_* \OO_{V^\prime}(-mE) \ = \ \overline{\aa^m},
\]
where as usual $\overline{ \aa^m}$ denotes the integral
closure of the ideal in question.  But quite generally, given any ideal
$\aa \subseteq \OO_V$, there exists an integer $c > 0$ such that
$\overline{\aa^{k+1}} = \aa \cdot \overline{\aa^k}$ for $k \ge c$ (cf.
\cite[(**) on p. 218]{PAG} and hence $\overline{\aa^m} \subseteq
\aa^{m-c}$ for $m \ge c$.
\end{proof}

\begin{rmk}
We do not know whether the normality hypothesis is essential.
\end{rmk}

\section{A   characterization of ample
divisors}\label{section:ampleness}

This section is devoted to  the statement and proof of our main
result.

We start by fixing  notation. In what follows,  $X$ will be a
projective variety of dimension $d$ over the complex numbers, and $L$
will denote a Cartier divisor on $X$.

If the divisor $L$ is  ample, then  the functions
$\^h^i$ vanish identically in a neighborhood of $[L]$ in $N^1(X)_\R$ for
every $i\geq 1$; this follows easily from  Serre vanishing, the
continuity of the functions $\^h^i$, and the fact that the ample cone is
open inside $N^1(X)_{\R}$. In particular these functions vanish at $L -
tA$ for every very ample divisor $A$ and every sufficiently small $t \ge
0$. The main result of this section is that this property characterizes
amplitude.

\begin{thm}\label{thm:ampleness} Let $X$ be a projective variety, and let
$L$ be a Cartier divisor on
$X$. Assume that there exists a very  ample divisor $A$ on
$X$ and a number $\ep > 0$ such that
$$
\^h^i(X,L-tA) \ = \ 0 \fall i > 0 \ , \ 0 \le t < \ep.
$$ Then $L$ is ample.
\end{thm}
\noi Theorem \ref{Max.Higher.Cohom} from the Introduction follows
immediately. We will deduce Corollary
\ref{Characterize.Ample.Classes.Cor.Intro} at the end of the section.

We now begin working towards the proof of Theorem \ref{thm:ampleness}.
First of all, in order to eventually be able to apply Proposition
\ref{lemma:base-loci}, we reduce  to the situation when the
variety $X$ is normal.

\begin{lem}\label{lemma:reduction-to-normal} Assume that
Theorem~\ref{thm:ampleness} holds for  normal projective varieties.
Then it  holds in general.
\end{lem}

\begin{proof} Let $X$ be an arbitrary projective variety, and suppose that
$L$ and $A$ are divisors on
$X$ satisfying the assumptions of Theorem~\ref{thm:ampleness}, so that
$L$ is Cartier, $A$ is ample, and there exists an $\ep > 0$ such that
$\^h^i(X,L-tA) = 0$ for all $i > 0$ and
$\ 0 \le t < \ep$. Consider the normalization $\n : \~X \to X$ of $X$.
Since
$\n$ is a finite morphism,
$\n^*A$ is ample. Moreover
\[\^h^i(\~X,\n^*(L-tA)) \ =\  \^h^i(X,L-tA). \]
 thanks to the birational invariance of  higher cohomology functions
(\cite[Proposition~2.9]{Kur}).
Assuming the theorem for normal varieties  we conclude that
$\n^*L$ is ample, and hence $L$ is ample as well.
\end{proof}

So we henceforth assume that $X$ is normal. The plan of the proof
is now to study the
$\widehat{h}^i$ via restrictions to divisors and use induction on
dimension. Specifically, choose a sequence of very general divisors
\[ E_1 \, ,\, E_2\, , \ldots \ \in \ \linser{A}.
\] Given $m, p > 0$ we take the first $p$ of the
$E_\alpha$ and form the complex $K^{\bullet}_{m,p}:$
\begin{equation}\label{Eqn:resolution}
\OO_X(mL) \lra \oversetoplus{p}\,
\OO_{E_\alpha}(mL) \lra
\oversetoplus{\binom{p}{2}} \, \OO_{E_\alpha
\cap E_\beta}(mL)
\lra \ldots\ ,
\end{equation} obtained as a twist of the $p$-fold tensor product of
the one-step complexes $\OO_X \lra \OO_{E_\alpha}$. Because it will  be
important   to keep track of the number of summands, we denote by
\[
\oversetoplus{\binom{p}{i}}
\OO_{E_{\alpha_1 }\cap
\ldots \cap E_{\alpha_i}}\]  the direct sum of the sheaves
$\OO_{E_{\alpha_1 }\cap
\ldots \cap E_{\alpha_i}}$ over all choices of $i$ increasing indices.  It
is established in
\cite[Corollary 4.2]{Kur} that $K^{\bullet}_{m,p}$ is acyclic, and hence
resolves
$\OO_X(mL - pA)$. In particular,
\begin{equation} \label{Cohom.is.Hypercohom.eqn}
\HH{r}{X}{\OO_X(mL - pA)} \ = \ {\mathbb{H}}^r
\big(  K^{\bullet}_{m,p} \big).
\end{equation}
The hypercohomology group on the right in
\eqref{Cohom.is.Hypercohom.eqn}   is in turn computed by a first-quadrant
spectral sequence with
\begin{equation} \label{Spect.Seq.Eqn}
 E_1^{i,j} \ = \
\begin{cases} H^j \Big( \OO_X(mL) \Big)  & i = 0 \\\oversetoplus{\binom{p}{i}}  H^j
\big(
\OO_{E_{\alpha_1 }\cap
\ldots \cap E_{\alpha_i}}(mL) \big) & i > 0.
\end{cases}
\end{equation}
As in \cite[2.2.37]{PAG} or \cite[Section 5]{Kur} we may
--- and do --- assume that the dimensions of all the groups appearing
on the right in \eqref{Spect.Seq.Eqn} are independent of the particular
divisors
$E_\alpha$ that occur. We will write these dimensions as
\[
\hhnospace{j}{\OO_{E_1}(mL)} \ \ , \ \
\hhnospace{j}{\OO_{E_1 \cap E_2}(mL)}
\] and so on.

The first point is to show that the
vanishing hypothesis of the theorem descends to very general divisors in
$\linser{A}$.
\begin{lem}\label{lemma:vanishing-on-E} Keeping  notation as in
Theorem~\ref{thm:ampleness}, assume that there is a positive real
number $\ep > 0$ such that
\begin{equation}
\label{Eqn.1}
\hhat{i}{X, L -tA} \ = \ 0 \fall i > 0 \ , \ 0
\le t < \ep.
\end{equation} Let $E \in \linser{A}$ be a very general divisor. Then
\begin{equation}\label{Eqn.2}
\^h^i(E, (L -tA)|_E) \ = \ 0 \fall i > 0 \ , \ 0 \le t <
\ep.
\end{equation}
\end{lem}

\begin{proof} Assuming \eqref{Eqn.1}, it's enough to prove
\begin{equation}
\label{Eqn.3}
\hhat{i}{E, L_E} \ = \ 0 \fall i > 0.
\end{equation} For then the more general statement
\eqref{Eqn.2} follows (using the homogeneity and continuity of the higher
cohomology functions on $X$ and on $E$) upon replacing $L$ by $L - \delta
A$ for a rational number $0 < \delta < \ep$.

Suppose then that \eqref{Eqn.3} fails, and consider the complex
$K^{\bullet}_{m,p}$. We compute a lower bound on the dimension of the
group
$E_{\infty}^{1,i}$ in the hypercohomology spectral sequence.
Specifically, by looking at the possible maps coming into and going out
from the $E_r^{1,i}$, one sees that
\begin{multline*}
\hhnospace{i+1}{X, mL - pA} \ + \
\hhnospace{i}{X, mL} \ \ge \\ p \cdot
\hhnospace{i}{\OO_{E_1}(mL)} \ - \
\binom{p}{2} \cdot \hhnospace{i}{\OO_{E_1 \cap E_2}(mL)} \ - \
\binom{p}{3} \cdot
\hhnospace{i-1}{\OO_{E_1 \cap E_2 \cap E_3}(mL)} \ - \
\ldots.
\end{multline*} Now we can find some fixed constant
$C_1 > 0$ such that for all $m \gg 0$:
\begin{align*}
\hhnospace{i}{\OO_{E_1 \cap E_2}(mL)} \ &\le  \ C_1
\cdot m^{d-2}, \\
\hhnospace{i-1}{\OO_{E_1 \cap E_2 \cap E_3}(mL)} \ &\le  \ C_1 \cdot
m^{d-3}, \
\text{etc.}
\end{align*} Moreover, since we are assuming for a contradiction that
$\hhat{i}{E, L_E} > 0$, we can find a constant $C_2 > 0$, together with a
sequence of arbitrarily large integers $m$, such that
\begin{equation}
\label{Eqn.4}
\hhnospace{i}{\OO_{E_1}(mL)} \ \ge \ C_2 \cdot m^{d-1}.
\end{equation} Putting this together, we find that there are arbitrarily
large integers $m$ such that
\begin{equation}\label{Eqn.5}
\hhnospace{i+1}{X,mL - pA} \, + \,
\hhnospace{i}{X, mL} \ \ge \ C_3 \cdot \Big( p m^{d-1} - p^2  m^{d-2} -
p^3   m^{d-3} -
\ldots \Big)
\end{equation} for suitable $C_3 >0$. Note that this constant $C_3$ is
independent of $p$. At this point, we fix a very  small rational  number
$0 <
\delta
\ll 1$. By the homogeneity of
$\^h^i$ on $E_1$, we can assume that the sequence of arbitrarily large
values of $m$ for which
\eqref{Eqn.4} and \eqref{Eqn.5} hold is taken among multiples of the
denominator of $\d$ (see  0.3). Then, restricting
$m$ to this sequence and taking $p =\d m$, the first term on the RHS of
\eqref{Eqn.5} dominates provided that $\delta$ is sufficiently small.
Hence
\[
\hhnospace{i+1}{X,mL - pA} \, + \,
\hhnospace{i}{X, mL} \ \ge \ C_4 \cdot \delta m^d
\] for a sequence of arbitrarily large $m$, and some
$C_4 > 0$. But this implies that
\[
\hhat{i+1}{X,L - \delta A} \  + \  \hhat{i}{X,L}
\ > \ 0,
\] contradicting the hypothesis.
\end{proof}

\begin{proof}[Proof of Theorem~\ref{thm:ampleness}] By
Lemma~\ref{lemma:reduction-to-normal}, we can suppose without loss of
generality that $X$ is normal. We
assume that there exists $\ep > 0$ such that
\begin{equation} \label{Eqn.6}
\hhat{i}{X, L - tA} \ = \ 0 \ \ \fall \ i > 0
\ , \ 0 \le t < \ep,
\end{equation} but that $L$ is not ample, and we'll aim to get a
contradiction.

Note that the Theorem fails for $L$ if and only if it fails for
integral multiples of $L - \delta A$ when $0 < \delta \ll 1$. So we can
suppose that we have a non-nef divisor $L$ satisfying \eqref{Eqn.6}.

Let $E \in \linser{A}$ be one of the very general divisors fixed at the
outset. Thanks to Lemma~\ref{lemma:vanishing-on-E}, we can assume by
induction on dimension that $\OO_E(L)$ is ample. Replacing $L$ again by a
multiple, we can suppose in addition that $\OO_E(L)$ is very ample with
vanishing higher cohomology. This combination being an open condition in
families, we can further assume that
$\OO_{E_\alpha}(L)$ is very ample for each of the
$E_\alpha$.

As above, form the complex
$K_{m,p}^{\bullet}$, and consider in particular the beginning of the
bottom row of the spectral sequence \eqref{Spect.Seq.Eqn}:
\[ \HHnospace{0}{\OO_X(mL)}
\overset{v_{m,p}}{\lra}
\oversetoplus{p}\,
\HHnospace{0}{\OO_{E_{\alpha}}(mL)}
\overset{u_{m,p}}{\lra}
\oversetoplus{\binom{p}{2}}\,
\HHnospace{0}{\OO_{E_\alpha \cap E_\beta}(mL)} .
\] There is a natural  injection
\[
\frac{\text{ker}(u_{m,p})}{\text{im}(v_{m,p})}
\ \subseteq \ \HH{1}{X}{mL - pA}, \]
and the plan is to estimate from below
the dimension of this subspace.

As in the proof of Lemma~\ref{lemma:vanishing-on-E} there is a uniform
bound having the shape
\[
\hhnospace{0}{\OO_{E_1 \cap E_2}(mL)} \ \le
\ C_1 \cdot m^{d-2}.
\] Therefore, considering $\text{ker}(u_{p,m})$ as a subspace of
$\oversetoplus{p}\,
\HHnospace{0}{\OO_{E_{\alpha}}(mL)}$, one has
\begin{equation} \label{Eqn.8}
\text{codim} \, \text{ker}(u_{m,p}) \ \le \ C_2 \cdot p^2 m^{d-2}
\end{equation} for some $C_2 > 0$ and all $m \gg 0$.

By Proposition~\ref{lemma:base-loci},   after possibly replacing $L$ by a
suitable multiple, we can find an ideal sheaf $\aa \subseteq \OO_X$
vanishing on a set of dimension $\ge 1$, together with  an integer
$c
\ge 0$ such that
$$
\bs{mL} \ \subseteq \ \aa^{m-c} \ \fall m > c.
$$ Then $v_{m,p}$ admits a factorization
\[
\xymatrix{
\HH{0}{X}{\OO_X(mL) \otimes \aa^{m-c}}
\ar[r]^(.46){v'_{m,p}} \ar@{=}[d]&
\oversetoplus{p}\,\HHnospace{0}{E_\a,\OO_{E_{\alpha}}(mL)
\otimes \aa^{m-c}}\ar@{^{(}->}[d] \\
\HH{0}{X}{\OO_X(mL)}
\ar[r]^(.46){v_{m,p}} & \oversetoplus{p}\,
\HHnospace{0}{E_\a,\OO_{E_{\alpha}}(mL)} }.
\] We claim that there is a constant $C_3 > 0$ such that for all $m \gg
0$:
\begin{equation} \label{Eqn.9}
\HHnospace{0}{\OO_{E_{\alpha}}(mL)
\otimes \aa ^{m-c}}\ \text{  has codimension } \ge C_3
\cdot m^{d-1}\ \text{ in }
\HHnospace{0}{\OO_{E_{\alpha}}(mL)}.
\end{equation} Granting this, we have
\begin{equation}\label{Eqn.10}
\dim
\frac{\text{ker}(u_{m,p})}{\text{im}(v_{m,p})}
\ \ge \ C_4 \cdot \Big( pm^{d-1} \ - \ p^2 m^{d-2} \Big)
\end{equation} for some constant $C_4 > 0$ and all $m
\gg 0$. Once again   fixing $0 < \delta \ll 1$, limiting $m$ to
multiples of the denominator of $\d$, and setting $p = \d m$, one finds
that
\[
\hh{1}{X}{\OO_X(m(L - \d A)) } \ \ge \ C_5 \cdot \delta m^d
\] for $m$ large enough. This implies that
$\hhat{1}{X, L - \delta A} > 0$, giving the required contradiction.

It remains to prove \eqref{Eqn.9}. Choose any point $x = x_{\alpha}
\in \text{Zeroes}(\aa) \cap E_\a$, and write $\mm_x
\subseteq \OO_{E_{\alpha}}$ for its maximal ideal: here we use that $\dim
\text{Zeroes}(\aa) \ge 1$ to know that such a point exists.   Since
\[ \HHnospace{0}{\OO_{E_{\alpha}}(mL)
\otimes \aa^{m-c}} \ \subseteq \
\HHnospace{0}{\OO_{E_{\alpha}}(mL)
\otimes \mm_x^{m-c}}, \] it is enough to bound the codimension of
$\HHnospace{0}{\OO_{E_{\alpha}}(mL) \otimes
\mm_x^{m-c}}$ in
$\HHnospace{0}{\OO_{E_{\alpha}}(mL)}$. But since
$\OO_{E_\a}(L)$ is very ample, it follows  that
$\OO_{E_\a}(mL)$ separates $(m-c)$-jets at the point
$x$. The dimension of the space of $(m-c)$-jets at a point of a possibly
singular variety is no smaller than the dimension of the space of
$(m-c)$-jets at a smooth point of a variety of the same
dimension, and thus we have
\[
\text{codim} \, \HHnospace{0}{\OO_{E_{\alpha}}(mL)
\otimes \mm_x^{m-c}} \ \ge \ \binom{m-c+d}{d-1},
\] as required.
\end{proof}

Finally, we give the proof of Corollary
\ref{Characterize.Ample.Classes.Cor.Intro} from the Introduction.
\begin{proof}[Proof of Corollary
~\ref{Characterize.Ample.Classes.Cor.Intro}] Consider the following three
statements concerning a Cartier divisor
$L$ on $X$:
\begin{enumerate}
\item
$L$ is ample;
\item for all $i > 0$ the function
$\^h^i$ vanishes in a neighborhood of $[L]$;
\item
$\hhat{i}{L-tA} = 0$ for some ample divisor $A$ and all
$i > 0$ and $0 \le t \ll 1$.
\end{enumerate} We have (a)~$\Rightarrow$~(b) by Serre's vanishing and
the continuity of
$\^h^i$, the implication (b)~$\Rightarrow$~(c) is obvious, and
Theorem~\ref{thm:ampleness} yields (c)~$\Rightarrow$~(a). Therefore
(a)~$\Leftrightarrow$~(b), which is the content of the
Corollary.
\end{proof}

\section{Amplitude of restrictions and cohomology}

In the course of this section we will study how the vanishing of higher
asymptotic cohomology relates to the amplitude of the restrictions of a
line bundle to certain very general complete intersections.

We consider as before an arbitrary complex projective variety $X$, and a
very ample divisor $A$ on $X$. For an arbitrary $\Q$-Cartier divisor $L$
on
$X$, we
 introduce the following invariants:
\begin{align*} a(L,A) \ &\deq \ \min \big \{ k \mid
\text{$L|_{E_1\cap\dots\cap E_k}$ is ample for very general $E_i \in
|A|$}\big \}, \\ b(L) \ &\deq \ \dim \mathbf{B}_+(L), \\ c(L) \ &\deq \ \max \big \{ i \mid
\text{$\^h^i$ is not identically zero in any neighborhood of $[L]$ in
$N^1(X)_\R$} \big \}.
\end{align*} These quantities express to some degree how far $L$ is from
being ample, with smaller values corresponding to 'more positive'
divisors. Note that they each depend only on the numerical equivalence
class of $L$.

The invariants $a(L,A),b(L),c(L)$ satisfy the following
relation.

\begin{prop}\label{proposition:c(L)<a(L,A)<b(L)} For every $L$ and $A$ as
above, we have
$$ c(L) \ \le \ a(L,A) \ \le \ b(L).
$$
\end{prop}

\begin{proof}[Proof of Proposition~\ref{proposition:c(L)<a(L,A)<b(L)}] The
inequality $a(L,A) \le b(L)$ follows from  the observation that the
restriction of a Cartier divisor to a general hyperplane section strictly
reduces the dimension of the augmented base locus, and the fact that a
divisor with empty augmented base locus is ample. Hence it remains to
show the inequality $c(L)
\le a(L,A)$.

At this point, it
will be  convenient to fix Cartier divisors $L_1,\dots,L_\r$ whose
numerical classes form a basis for $N^1(X)_\R$, and to work in the
computations that follow with
$\Q$-divisors on $X$ are taken to  be linear combinations of the $L_k$.
The point of this is that  given    very general divisors
\[ E_1 \ ,\ E_2\ , \ \ldots \ \in \linser{A},
\] we can assume as above that for  every $\vec{m} =(m_k) \in \Z^\r$, the
dimensions of the cohomologies  of the restriction of $\sum m_kL_k$ to any
intersection of
$j$ distinct
$E_\a$ are independent of the particular $E_\a$ chosen. We will
implicitly make use of this assumption in the following paragraphs.

We now begin the proof that  $c(L)
\le a(L,A)$. By replacing $L$ by an integral multiple of
an arbitrarily close rational class in $N^1(X)_\R$, it suffices to show
that
$\hhat{i}{X,L} = 0$ for every $i > a(X,L)$. The plan is to proceed by
induction on
$a(L,A)$. If $a(L,A) = 0$, then the claim follows by Serre's vanishing, so we
can assume that $a(L,A) \ge 1$ and suppose that the property holds for
$a(L,A)-1$.

We claim that for every $i > a(L,A)$ the function
$$ g_i(t) \ \deq \ \hhat{i}{L + tA}
$$ is differentiable with zero derivative on
$(0,\infty)$. Granting this for the moment, we will finish the proof. Indeed
we have $g_i(t) = 0$ for $t \gg 0$ by Serre's vanishing, and thus the
claim implies that
$g_i$ vanishes on $(0,\infty)$, hence at $0$ by continuity.

Turning to the proof of the claim, consider to begin with an arbitrary
divisor $D$ on $X$, and arbitrary integers $m , p \ge 0$. By looking again at the
spectral sequence \eqref{Spect.Seq.Eqn}, one sees that
\begin{equation} \label{diff.quot.estimate.eqn}
\begin{aligned}
\Big|\,\hhnospace{i}{mD-pA} -
\hhnospace{i}{mD}\,\Big| \ & \le \ p\.
\big (\hhnospace{i}{mD|_{E_1}} +
\hhnospace{i-1}{mD|_{E_1}}\big ) \\ + \ \binom{p}{2}\. &
\big (\hhnospace{i-1}{mD|_{E_1 \cap E_2}} +
\hhnospace{i-2}{mD|_{E_1 \cap E_2}}\big ) + \dots .
\end{aligned}
\end{equation}
Now fix rational numbers $t, \delta >0$, with $ \delta \ll1$, and set $p
=
\delta m$ where
$m$ clears the demoninators of $t$ and $\delta$. Applying
\eqref{diff.quot.estimate.eqn} with
$D = L + tA$, and taking $m$ to range
over a sequence of suitably divisible integers computing the
largest of
$\hhat{i}{D}$ and $\hhat{i}{D-\d A}$, one finds  that
\begin{equation} \label{Eqn:12}
\begin{aligned}
\frac{\big|\, \hhat{i}{D-\d A} - \hhat{i}{D}\, \big|}{\d} \ & \le \
C_1\.\(\hhat{i}{D|_{E_1}} + \hhat{i-1}{D|_{E_1}}\) \\ +
\ C_2\.\d\.&\(\hhat{i-1}{D|_{E_1 \cap E_2}} +
\hhat{i-2}{D|_{E_1 \cap E_2}}\) + \dots ,
\end{aligned}
\end{equation}
 where the $C_j$ are  positive numerical constants
depending only  on the dimension
$d$.  A similar argument
using $D  + \d A$ in place of $D$ gives
\begin{align}\label{Eqn:13}
\frac{\big|\, \hhat{i}{D+\d A)} - \hhat{i}{D}\, \big |}{\d} \ & \le \
C_1\.\(\hhat{i}{(D+\d A)|_{E_1}} + \hhat{i-1}{(D+\d A)|_{E_1}}\) \\ + \
C_2\.\d\. &\(\hhat{i-1}{(D+\d A)|_{E_1 \cap E_2}} + \hhat{i-2}{(D+\d
A)|_{E_1 \cap E_2}}\) + \dots . \nonumber
\end{align}

Suppose now that $i > a(L,A)$, and as above let $D = L + tA$.  Since
$$ a(D|_{E_1 \cap\dots\cap E_k}) \ = \ a(D,A) - k
\ \le \ a(L,A) - k \ < \ i - k
$$ for all $1 \le k \le a(D,A)$, the same holding with
$D+\d A$ in place of $D$, all the asymptotic cohomological dimensions
other than $\^h^0$ appearing in the RHS of \eqref{Eqn:12} and
\eqref{Eqn:13} are zero by induction, and hence each of the dimensions
$\^h^0$ compute the self-intersection of the given divisor. Thus we can
find a positive constant $C$, independent of $\d$, such that
\begin{equation} \label{Diff.Quot.Estimate.eqn}
\frac{\big| \,\hhat{i}{L \pm \d A} - \hhat{i}{L}\, \big|}{\d} \ \le \
C\.\d^{i-1}\.\((L+\d A)^{d-i}\.A^i + \d\.(L+\d A)^{d-i-1}\.A^{i+1}\) .
\end{equation}
So far we have assumed that $\delta \in \Q$, but then  by the continuity
of both sides,  the inequality must  actually hold for  every real $0 < \d
\ll 1$. Therefore $g_i(t)$ is differentiable with zero derivative at
every rational $t >0$. To see that this property extends over the whole
interval
$(0,\infty)$, we fix an arbitrary $t > 0$, and let
$\d$ be any real number with $0 < |\d| \ll 1$ and such that $t + \d \in
\Q$. Then \eqref{Diff.Quot.Estimate.eqn} yields
\begin{multline*}
\frac{\Big{|} \, \hhat{i}{(L + t A) + \d A} -
\hhat{i}{L+tA}\, \Big{|}}{| \d |}
\ = \ \frac{\Big{|} \, \hhat{i}{(L + (t+\d) A) - \d A} -
\hhat{i}{L+(t + \d) A\,} \Big{|}}{| \d |} \\
\le \ C\.|\d|^{i-1}\.\((L+(t+\d+|\d|)A)^{d-i}\.A^i +
|\d|\.(L+(t+\d+|\d|)A)^{d-i-1}\.A^{i+1}\).
\end{multline*}   Since the RHS goes to zero as $\d \to 0$, this implies
that
$g_i(t)$ is differentiable at $t > 0$ and, moreover, we have
$g_i'(t) = 0$. This proves the claim, hence completes the proof of the
theorem.
\end{proof}

One might ask whether in fact equality holds in $c(L)\leq a(L,A)$ for all very ample
divisors $A$. Although the idea is tempting, it turns out that in general
this is not the case.

\begin{eg}\label{eg:c(L)<a(L,A)<b(L)} Let $X = \F_1
\times \P^1$, where $\F_1 \cong \Bl_p \P^2$, and denote by $p : X \to
\F_1$ and $q : X \to \P^1$ the two projections. Let $E \subset \F_1$ be
the $(-1)$-curve, let $F \subset \F_1$ be a fiber of the ruling, and let
$H \subset \P^1$ be a point. We consider the divisors
$$ L = p^*(\l E + F) + q^*H \ \andd \ A = p^*(E + \m F) + q^*H \quad
\text{for some} \quad \l, \m \in\Z_{\ge 2}.
$$ Notice that $A$ is very ample and $L$ is a big divisor. Moreover, the
base locus of $L$ coincides with its augmented base locus and is equal to
$B \deq p^{-1}(E)$. In particular $b(L) = 2$.

On the other hand,
K\"unneth's formula for asymptotic cohomological functions (see
\cite[Remark~2.14]{Kur}), and the fact that $L$ is not ample imply that
$c(L) = 1$. Fix a general element $Y \in \linser{A}$ cutting out a smooth
divisor $D$ on $B$. Note that
$\O_B(D) \cong \O_{\P^1\times\P^1}(\m-1,1)$ via the isomorphism $B = E
\times \P^1 \cong \P^1 \times \P^1$. Therefore, since $D$ is smooth, it
must be irreducible; moreover,
$p$ induces an isomorphism $D \cong \P^1$.
We observe that the base locus
of $L|_Y$ is contained in the restriction of the base locus of $L$, hence
in $D$, and that $\O_D(L|_Y) \cong \O_{\P^1}(\m - \l)$. We conclude that
$$ a(L,A) \ = \
\begin{cases} 2 = b(L) &\text{if $\l \ge \m$}, \\ 1 = c(L) &\text{if $\l
< \m$}.
\end{cases}
$$
\end{eg}

\begin{rmk} Outside the big cone, the equation $a(L,A) = c(L)$ fails
badly on every projective variety. Indeed, if $L$ is a Cartier divisor on
a $d$-dimensional projective variety $X$ such that neither $L$ nor $-L$
is pseudo-effective, then we have $c(L) < d$ and $c(-L) < d$, and either
$a(L) = d$ or $a(-L) = d$.
\end{rmk}

\end{document}